\documentclass[12pt,reqno]{amsart}
\usepackage{latexsym,color,amssymb,eucal,epsfig,bbm,amsfonts} 

\begin{document} 

\newenvironment{mymatrix}{\begin{array}{*{20}{c}}}{\end{array}}

\renewcommand{\qed}{$\blacksquare$} 

\makeatletter
\renewcommand*\env@matrix[1][c]{\hskip -\arraycolsep
  \let\@ifnextchar\new@ifnextchar
  \array{*\c@MaxMatrixCols #1}}
\makeatother

\newcommand{\F}{\mathbb{F}}
\newcommand{\Q}{\mathbb{Q}}
\newcommand{\R}{\mathbb{R}}
\newcommand{\C}{\mathbb{C}}
\newcommand{\N}{\mathbb{N}}
\newcommand{\Z}{\mathbb{Z}}

\newcommand{\gln}{{\mathbf{GL}_n}}
\newcommand{\Log}{\mathrm{Log}} 
\newcommand{\np}{\mathbf{NP}} 
\newcommand{\cO}{{\mathcal{O}}}
\newcommand{\area}{{\mathrm{Area}}}
\newcommand{\anewt}{{\mathrm{ArchNewt}}}
\newcommand{\nanewt}{{\widetilde{\anewt}}}
\newcommand{\viro}{{\mathrm{Viro}}} 
\newcommand{\wiro}{{\widetilde{\viro}}}
\newcommand{\sign}{\mathrm{sign}}
\newcommand{\belta}{\overline{\Delta}}
\newcommand{\thth}{^{\text{\underline{th}}}}
\newcommand{\eps}{\varepsilon}
\newcommand{\var}{\mathrm{Var}}
\newcommand{\amoeba}{\mathrm{Amoeba}}
\newcommand{\dia}{$\diamond$}
\newcommand{\Qn}{\Q^n}
\newcommand{\cC}{\mathcal{C}}
\newcommand{\cW}{\mathcal{W}}
\newcommand{\bC}{\bar{\cC}}
\newcommand{\cI}{\mathcal{I}}
\newcommand{\cJ}{\mathcal{J}}
\newcommand{\cR}{\mathcal{R}}
\newcommand{\Rn}{\R^n}
\newcommand{\Rm}{\R^m}
\newcommand{\Cn}{\C^n}
\newcommand{\bO}{\mathbf{O}}
\newcommand{\cA}{{\mathcal{A}}}
\newcommand{\hA}{\hat{\cA}}
\newcommand{\cS}{\mathcal{S}}
\newcommand{\Pro}{{\mathbb{P}}}
\newcommand{\barf}{\bar{f}}
\newcommand{\babla}{\overline{\nabla}}
\newcommand{\Zn}{\Z^n}
\newcommand{\Cm}{\C^m}
\newcommand{\Cs}{\C^*}
\newcommand{\Csn}{(\Cs)^n}
\newcommand{\Rs}{\R^*}
\newcommand{\cF}{{\mathcal{F}}}
\newcommand{\cV}{{\mathcal{V}}}
\newcommand{\conv}{{\mathrm{Conv}}}
\newcommand{\archnewt}{{\mathrm{ArchNewt}}}
\newcommand{\supp}{{\mathrm{Supp}}}

\newtheorem{thm}{Theorem}[section]
\newtheorem{algor}[thm]{Algorithm} 
\newtheorem{cor}[thm]{Corollary} 
\newtheorem{lemma}[thm]{Lemma}
\newtheorem{prop}[thm]{Proposition}
\newtheorem{question}[thm]{Question}
\newtheorem{conj}[thm]{Conjecture}
\newtheorem{problem}[thm]{Problem}
\theoremstyle{definition}
\newtheorem{dfn}[thm]{Definition}
\newtheorem{rem}[thm]{Remark}
\newtheorem{ex}[thm]{Example}
\newtheorem{prt}{Passare-Rullg\aa{}rd Theorem}
\renewcommand{\theprt}{\unskip}
\newtheorem{hk}{The Horn-Kapranov Uniformization}
\renewcommand{\thehk}{\unskip}
\newtheorem{kapnona}{Archimedean Amoeba Theorem}
\renewcommand{\thekapnona}{\unskip}
\theoremstyle{remark}
\newtheorem{remark}[thm]{Remark}

\numberwithin{equation}{section}

\title[Randomization and Faster Real Root Counting]{
\mbox{}\\
\vspace{-1in}
Randomization, Sums of Squares, and Faster 
Real Root Counting for Tetranomials and Beyond\\ \today } 

\author[O.\ Bastani]{Osbert Bastani}
\author[C.\ Hillar]{Christopher J.\ Hillar}
\author[D.\ Popov]{Dimitar Popov}
\author[J.\ M.\ Rojas]{J.\ Maurice Rojas}
\address{Harvard University, Massachusetts Hall, Cambridge, MA \ 02138} 
\email{hypo3400@gmail.com}
\address{Redwood Center for Theoretical Neuroscience, 575A Evans Hall, MC 3198
Berkeley, CA \ 94720-3198. }
\email{chillar@msri.org}
\address{MIT, 77 Mass.\ Ave., Cambridge, MA \ 02139} 
\email{dpopov@mit.edu} 
\address{TAMU 3368, Texas A\&M University, College Station, TX \ 77843-3368 } 
\email{rojas@math.tamu.edu} 
\thanks{Bastani and Popov were partially supported by NSF REU grant 
DMS-0552610. Hillar was partially supported by an NSF Postdoctoral 
Fellowship and an NSA Young Investigator grant. Rojas was partially supported 
by NSF CAREER Grant DMS-0349309, 
a Wenner Gren Foundation grant, Sandia National Laboratories, and MSRI}

\begin{abstract} 
Suppose $f$ is a real univariate polynomial of degree $D$ with 
exactly $4$ monomial terms. We present an algorithm, with 
complexity polynomial in $\log D$ on average (relative to the stable 
log-uniform measure), for counting the number of real roots of $f$.  
The best previous algorithms had complexity super-linear in $D$. 
We also discuss connections to sums of squares and $\cA$-discriminants, 
including explicit obstructions to expressing 
positive definite sparse polynomials as sums of squares of few sparse 
polynomials. Our key tool is the introduction of efficiently computable 
{\em chamber cones}, bounding regions in coefficient space where the
number of real roots of $f$ can be computed easily. Much of our theory 
extends to $n$-variate $(n+3)$-nomials.
\end{abstract} 
\maketitle

\vspace{-.7cm}
\section{Introduction} 
Counting the real solutions of polynomial equations in one variable 
is a fundamental ingredient behind many deeper tasks and applications 
involving the topology of real algebraic sets. However, the intrinsic 
complexity of this basic enumerative problem becomes a mystery as soon 
as one considers the input representation in a refined way. 
Such complexity questions have practical impact for, in many  
applications such as geometric modelling or the discretization of 
physically motivated partial differential equations, one encounters 
polynomials that have sparse expansions relative to some basis. 
So we focus on new, exponential speed-ups for counting the real roots 
of certain sparse univariate polynomials of high degree. 

Sturm sequences \cite{sturm}, and their later refinements 
\cite{habicht,bpr}, have long been a centrally important technique for 
counting real roots of univariate polynomials. In combination with more 
advanced algebraic tools such as a Gr\"obner bases or resultants 
\cite{gkz94,bpr}, Sturm sequences have even been applied to algorithmically 
study the topology of real algebraic sets in arbitrary dimension 
(see, e.g., \cite[Chapters 2, 5, 11, and 16]{bpr}). However, as 
we will see below (cf.\ Examples \ref{ex:tri} and \ref{ex:tetra1}), there 
are obstructions to attaining polynomial intrinsic complexity,  
for {\em sparse} polynomials, via Sturm sequences. So we must seek 
alternatives. 

More recently, relating multivariate positive polynomials to sums of squares 
has become an important algorithmic tool in optimizing real polynomials over 
semi-algebraic domains \cite{parrilo,lasserre2}. However, there are 
also obstructions to the use of sums of squares toward speed-ups for 
sparse polynomials (see Theorem \ref{thm:neg} below). 

Discriminants have 
a history nearly as long as that of Sturm sequences and sums of squares, but 
their algorithmic power has not yet been fully exploited. 
Our main result is that {\em $\cA$-discriminants} \cite{gkz94} yield an 
algorithm for counting real roots, with average-case 
complexity polynomial in the {\em logarithm} of the degree,  
for certain choices of probability distributions on the input  
(see Theorem \ref{thm:pos} below). 
The use of randomization is potentially inevitable in light of the 
fact that even detecting real roots becomes $\np$-hard already for 
moderately sparse multivariate polynomials \cite{brs,prt}.    

\subsection{From Large Sturm Sequences to Fast Probabilistic Counting} 
The classical technique of Sturm Sequences \cite{sturm,bpr} reduces 
counting the roots of a polynomial $f$ in a half-open interval $[a,b)$ 
to a gcd-like computation, followed by sign evaluations for a 
sequence of polynomials. A key difficulty in these methods, however, 
is their apparent super-linear dependence on the degree of the 
underlying polynomial. 
Consider the following two examples (see also \cite[Example 1]{rojasye}). 
\begin{ex} 
\label{ex:tri} {\em 
Setting $f(x_1)\!:=\!x^{317811}_1-2x^{196418}_1+1$, the {\tt realroot} 
command in {\tt Maple 14}\footnote{Running on a 16GB RAM Dell PowerEdge 
SC1435 departmental server with 2 dual-core Opteron 2212HE 2Ghz processors 
and OpenSUSE 10.3.}  (which is an implementation of Sturm 
Sequences) results in an out of memory error after about $31$ seconds. 
The polynomials in the underlying computation, while quite sparse, have 
coefficients with hundreds of thousands of digits, 
thus causing this failure. On the other hand, 
via more recent work \cite{brs}, one can show that when $c\!>\!0$ 
and $g(x_1)\!:=\!x^{317811}_1-cx^{196418}_1+1$, $g$ has exactly $0$, $1$, or 
$2$ positive roots according as $c$ is less than, equal to, or greater than 
$\frac{317811}{(121393^{121393} 
196418^{196418})^{1/317811}}\!\approx\!1.944...$. 
In particular, our $f$ has exactly $2$ positive roots. (We discuss how to 
efficiently decide the size of monomials in rational numbers 
with rational exponents in Algorithm \ref{algor:sign} of Section 
\ref{sub:bit} below.) \dia } 
\end{ex} 
\begin{ex} 
\label{ex:tetra1} 
{\em 
Going to tetranomials, consider  
$f(x_1)\!:=\!ax^{100008}_1-x^{50005}_1+bx^{50004}_1-1$ with $a,b\!>\!0$. 
Then (via the classical Descartes' Rule of Signs \cite[Cor.\ 10.1.10, pg.\ 319]
{ana}) such an $f$ has exactly $1$ or $3$ positive roots, but the 
inequalities characterizing which $(a,b)$ yield either possibility are 
much more unwieldy than in our last example: there are at least $2$, involving  
polynomials in $a$ and $b$ having tens of thousands of terms. 
In particular, for $(a,b)\!=\!\left(2,\frac{1}{2}\right)$, Sturm sequences on 
{\tt Maple 14} result in an out of memory error after about 122 seconds. 
\dia } 
\end{ex}

We have discovered that $\cA$-discriminants, reviewed in 
Section \ref{sec:back}, enable algorithms with complexity 
polynomial in the {\em logarithm} of the degree. 
\begin{dfn}
{\em 
For any $x\!=\!(x_1,\ldots,x_d)\!\in\!\C^d$, we define 
$\Log|x|\!:=\!(\log|x_1|,\ldots, \log|x_d|)$ and let the 
{\em stable log-uniform measure} on $\R^d_+$ (resp.\ $\Z^d$) be the 
probability measure $\nu$ (resp.\ $\nu'$) defined 
as follows: 
$\nu(S)\!:=\!\lim\limits_{M\rightarrow +\infty} \frac{\mu(\Log|S|\cap 
[-M,M]^d)}{(2M)^d}$ (resp.\ $\nu'(S)\!:=\!\lim\limits_{M\rightarrow +\infty} 
\frac{\#(\Log|S|\cap \{-M,\ldots,M\}^d)}{(2M+1)^d}$), where 
$\mu$ denotes the standard Lebesque measure on $\R^d$ and 
$\#(\cdot)$ denotes set cardinality. \dia}  
\end{dfn} 

\noindent 
Note that the stable log-uniform measure is finitely additive (but not 
countably additive), and is invariant under reflection across coordinate 
hyperplanes.
\begin{thm} 
\label{thm:pos} 
Let $0\!<\!a_2\!<\!a_3\!<\!a_4\!=\!D$ be positive integers,\\ 
\mbox{}\hfill 
$f(x_1)\!:=\!c_1+c_2 x^{a_2}_1+c_3 x^{a_3}_1+c_4 x^{a_4}_1$ 
\hfill \mbox{}\\  
with $c_1,\ldots,c_4$ being independent stable log-uniform  
random variables chosen from $\R$ (resp.\ $\Z$), and 
define $h\!:=\!\log(2+\max_i|c_i|)$. 
Then there is a deterministic algorithm, with  
complexity polynomial in $\log D$ (resp.\ 
$h$ and $\log D$), that computes a number in $\{0,1,2,3\}$ that, 
with probability $1$, is exactly the number of real roots of $f$.  
The underlying computational model is the BSS model over $\R$ (resp.\ the 
Turing model).  
\end{thm}  

\noindent 
The key idea is that while the regions of coefficients determining 
polynomials with a constant number of real roots become more complicated as 
the number of monomial terms increases, one can nevertheless efficiently 
characterize large subregions --- {\em chamber cones} --- where the number of 
real roots is very easy to compute.  This motivates the introduction of 
probability and 
average-case complexity.  The $\cA$-discriminant allows one to make this 
approach completely precise and algorithmic. In fact, our framework enables 
us to transparently extend Theorem \ref{thm:pos} to $n$-variate $(n+3)$-nomials 
(see Theorem \ref{thm:nplus3} of Section \ref{sub:can}). 

Our focus on the stable log-uniform measure simplifies our development 
and has some practical motivation: when one considers $N$-bit floating-point 
numbers with uniformly random exponent and mantissa, taking 
$N\longrightarrow +\infty$ and suitably rescaling 
yields exactly the stable log-uniform measure on $\Z$. The stable 
log-uniform measure has also been used in work of Avenda\~{n}o and 
Ibrahim to study the expected number of roots of sparse polynomial systems 
over local fields other than $\R$ \cite{ave}. 

It is of course quite natural to ask how the expected complexity in 
Theorem \ref{thm:pos} behaves under 
other well-known measures, e.g., uniform or Gaussian. Unfortunately, the 
underlying calculations become much more complicated. On a deeper level, it 
is far from clear what a truly ``natural'' probability measure on the space of 
tetranomials is. For instance, for non-sparse polynomials, it is 
popular to use specially weighted independent Gaussian coefficients since 
the resulting measure becomes invariant under a natural orthogonal 
group action (see, e.g., \cite{kostlan,ss4,bsz}). 
However, we are unaware of any study of the types of distributions  
occuring for the coefficients of polynomials actually occuring in physical 
applications. 

The speed-ups we derive actually hold in far greater generality:  
see \cite{brs,prt} for the case of $n$-variate $(n+k)$-nomials 
with $k\!\leq\!2$, Section \ref{sec:chambercones} here for 
connections to $n$-variate $(n+3)$-nomials, the forthcoming paper \cite{aar} 
for the general univariate case, and the forthcoming paper \cite{prrt} for 
chamber cone theory for $n\times n$ sparse polynomial systems. 
One of the main goals of our paper is thus to illustrate 
and clarify the underlying theory in a non-trivial special case. We now 
state our second main result. 

\subsection{Sparsity and Univariate Sums of Squares} 
Recent advances in semidefinite programming have produced efficient
algorithms for finding sum of squares representations of certain nonnegative
polynomials, thus enabling efficient polynomial optimization under 
certain conditions. When the input is a sparse polynomial it is then natural 
to ask if there is a sum of squares representation that also respects 
sparseness. Indeed, it is well-known that a nonnegative univariate 
polynomial can always be written as a sum of two squares of, usually 
non-sparse, polynomials (see, e.g., \cite{pourchet} for 
refinements). The following result 
demonstrates that a sparse analogue is either unlikely or much more subtle. 
\begin{thm} 
\label{thm:neg} 
There do {\em not} exist absolute constants $\ell$ and $m$ with the 
following property: Any trinomial $f\!\in\!\R[x_1]$ 
that is positive on $\R$ can be written in the form 
$f\!=\!g^2_1+\cdots+g^2_\ell$, for some $g_1,\ldots,g_\ell\!\in\!\R[x_1]$ 
with $g_i$ having at most $m$ terms for all $i$. 
\end{thm} 

\noindent 
Were there to be a sufficiently efficient representation of positive sparse 
polynomials as sums of squares, one could then try to use semidefinite 
programming to find such a representation explicitly for a given 
polynomial. This in turn could yield an efficient reduction from deciding the 
existence of real roots to a (small) semidefinite programming problem, 
similar to the techniques of \cite{parrilo}. Our last theorem thus reveals an 
obstruction to this sums of squares approach. 

\subsection{Related Approaches} 
The best known algorithms for real root 
counting lack speed-ups for sparse polynomials like the 
average-case complexity bound from our first main result. For example, 
in the notation of Theorem \ref{thm:pos}, 
\cite{lickroy} gives an arithmetic complexity bound of 
$O(D\log^5 D)$ which, via the techniques of \cite{bpr}, yields 
a bit complexity bound super-linear in $h+D$. No algorithm with 
complexity polynomial in $\log D$ (deterministic, randomized, or 
average-case) appears to have been known before for tetranomials. 
(See \cite{htzekm} for recent speed benchmarks of univariate 
real solvers.) 

As for alternative approaches, softening our concept of sparse sum of 
squares representation may still enable speed-ups similar to Theorem  
\ref{thm:pos} via semidefinite programming. For instance, one could ask if a 
positive  trinomial of degree $D$ always admits a representation as a sum 
of $\log^{O(1)}\!D$ squares of polynomials with $\log^{O(1)}\!D$ terms. 
This question appears to be completely open. 
\begin{ex}\label{ex:log} {\em Observe that a quick derivative computation 
yields that\\
\mbox{}\hfill   
$f(x_1)\!:=\!x^{2^k}_1-2^kx_1+2^{k}-1$\hfill\mbox{}\\ 
attains a unique minimum value of 
$0$ at $x\!=\!1$. So this $f$ is nonnegative. On the other hand,  
one can prove easily by induction that $f(x_1) = 2^{k-1}
\sum_{i=0}^{k-1} \frac{1}{2^i}\left(x^{2^i}_1-1\right)^2$, 
thus yielding an expression for $f$ as a sum of $O(\log D)$ {\em bi}nomials 
with $D\!=\!2^k$. \dia } 
\end{ex} 

Note also that while we focus on speed-ups that replace the polynomial 
degree $D$ by $\log D$ in this paper, other practically important speed-ups 
combining semidefinite programming and sparsity are certainly possible 
(see, e.g., \cite{lasserre1,kojima}).   

\section{Background}
\label{sec:back} 

\subsection{Amoebae and Efficient $\cA$-Discriminant Parametrization} 
Let us first briefly\linebreak review two important constructions by 
Gelfand, Kapranov, and Zelevinsky. 
\begin{dfn}
\label{dfn:first} 
{\em Let $c\!:=\!(c_1,\ldots,c_m)$, let 
$\cA\!=\!\{a_1,\ldots,a_m\}\!\subset\!\Zn$ 
have cardinality $m$, and define the corresponding family of 
(Laurent) polynomials\\
\mbox{} \hfill $\cF_\cA\!:=\!\{ c_1x^{a_1}+\cdots+c_mx^{a_m} \; | \; 
c\!\in\!\C^m\}$,\hfill \mbox{}\\
where the notation $x^{a_i}\!:=\!x^{a_{1,i}}_1\cdots x^{a_{n,i}}_n$ 
is understood. When $c_i\!\neq\!0$ for all
$i\!\in\!\{1,\ldots,m\}$ then we call $\cA$ the {\em support} of
$f(x)\!=\!\sum^m_{i=1} c_ix^{a_i}$, also using the notation
$\supp(f)$. \dia} 
\end{dfn}

\noindent 
\begin{minipage}[b]{0.7\linewidth}
\vspace{0pt}
\begin{dfn} 
\label{dfn:amoeba} 
{\em 
For any field $K$ we let $K^*\!:=\!K\setminus\{0\}$. 
Given any $g\!\in\!\C[x_1,\ldots,x_n]$, we then  
define its {\em amoeba}, $\amoeba(g)$, to be $\{\Log|c| \; \ | \ \; 
c\!=\!(c_1,\ldots,c_m)\!\in\!(\C^*)^m \text{ and } 
g(c_1,\ldots,c_m)\!=\!0\}$. \dia} 
\end{dfn} 
\begin{kapnona} 
(weaker version of \cite[Cor.\ 1.8]{gkz94}) 
{\em Following the notation of Definition \ref{dfn:amoeba}, 
the complement of $\amoeba(g)$ in $\R^m$ is  
a finite disjoint union of open convex sets. \qed} 
\end{kapnona}  

\noindent 
An example of an amoeba of a bivariate polynomial (see 
Example \ref{ex:tetra2} below) appears in the right-hand 
illustration. While the complement of the amoeba (in white) appears to 
have $3$ convex connected components, there are in fact $4$: the fourth 
component is a thin sliver emerging further below from the downward pointing 
tentacle. 
\end{minipage} 
\begin{minipage}[b]{0.3\linewidth}
\vspace{0pt}
\mbox{}\hfill\raisebox{0cm}{\epsfig{file=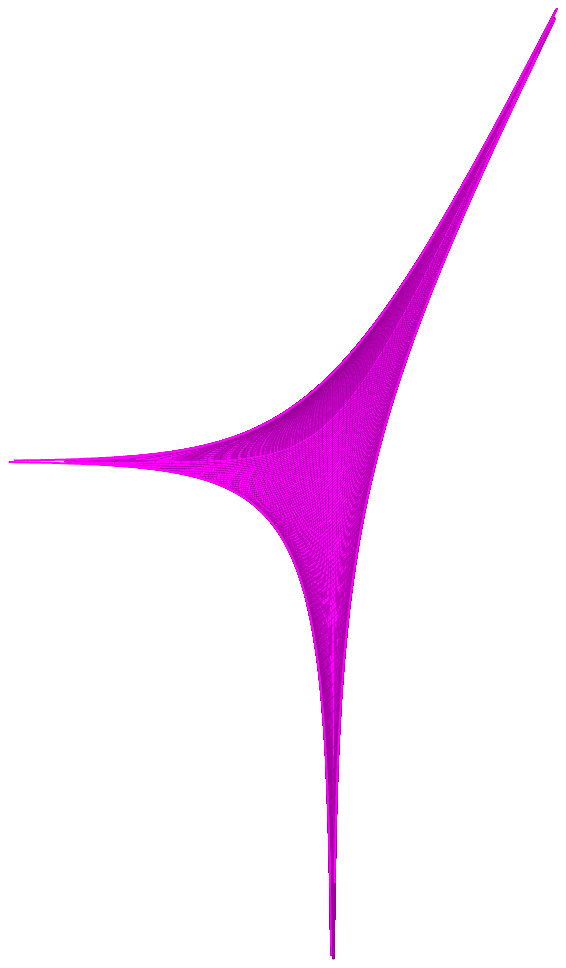,
height=2.6in,clip=}} 
\end{minipage} 
\begin{dfn}
{\em Following the notation of Definition \ref{dfn:first} and 
letting $f(x)\!:=\!c_1x^{a_1}+\cdots +c_mx^{a_m}$, 
we define $\nabla_\cA$ --- the
{\em $\cA$-discriminant variety} \cite[Chs.\ 1 \& 9--11]{gkz94} --- to be the
closure of the set of all $[c_1:\cdots :c_m]\!\in\!\Pro^{m-1}_\C$
such that\\ 
\mbox{}\hfill $f=\frac{\partial f}{\partial x_1}=
\cdots=\frac{\partial f}{\partial x_n}=0$\hfill\mbox{}\\ 
has a solution in $\Csn$. We then define (up to sign) the {\em 
$\cA$-discriminant}, $\Delta_\cA\!\in\!\Z[c_1,\ldots,c_m]$, to be the
(irreducible) defining polynomial of $\nabla_\cA$ when $\nabla_\cA$ is a 
hypersurface. Finally, we 
let $\nabla_\cA(\R)$ denote the real part of $\nabla_\cA$. \dia } 
\end{dfn}
\begin{rem} \label{rem:gen} {\em The $\nabla_\cA$ considered in this paper 
will all ultimately be hypersurfaces. \dia} 
\end{rem} 
\begin{ex}
\label{ex:tetra2}
{\em 
Taking $\cA\!=\!\{0,404,405,808\}$, we see that\\
\begin{minipage}[b]{0.62\linewidth}
$\cF_\cA$ consists simply of polynomials of the form
$f(x_1)\!:=$\linebreak
$c_1+c_2x^{404}_1+c_3x^{405}_1+c_4x^{808}_1$.
The underlying $\cA$-discriminant is then a polynomial in the $c_i$ having
$609$ monomial terms and degree $1604$.
However, while $\Delta_\cA$ is unwieldy, we can still 
easily plot the real part of its zero set $\nabla_\cA(\R)$ via the 
{\em Horn-Kapranov\linebreak 
Uniformization} (see its statement below, and the \linebreak
illustration to the right). \dia   
\end{minipage}
\hspace{.5cm}
\begin{picture}(100,10)(0,0)
\put(15,-3){\epsfig{file=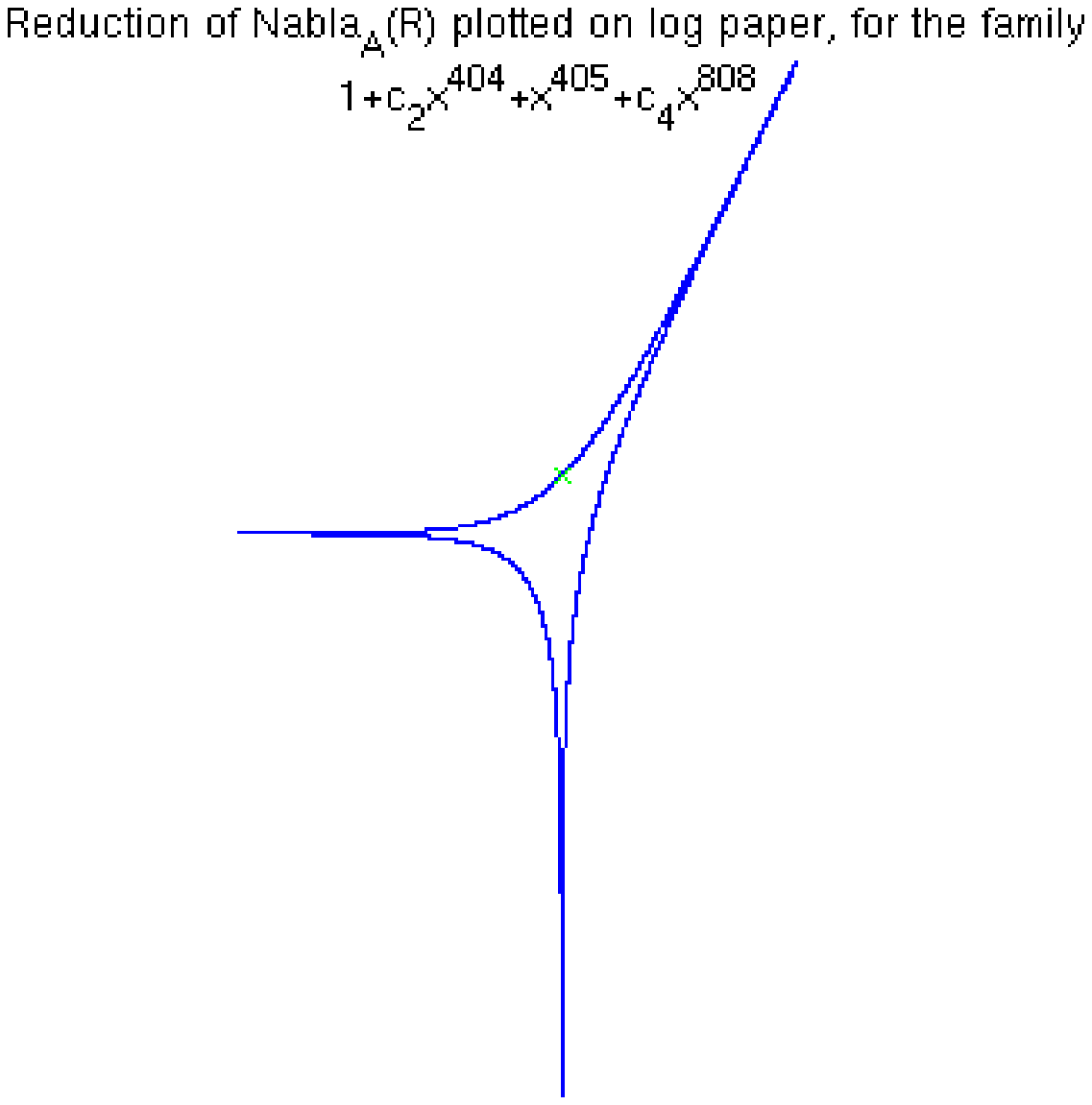,height=1.7in}}  
\end{picture}} 
\end{ex}

\noindent 
The plotted curve above is the image of the real roots of
$\belta_\cA(c_2,c_4)\!:=\!\Delta_\cA(1,c_2,1,c_4)$ under 
the $\Log|\cdot|$ map, i.e., the amoeba of $\belta_\cA$. 
Amoebae give us a convenient way to introduce polyhedral/tropical methods 
into our setting. For our last example, the boundary of $\amoeba(\belta_\cA)$ 
is contained in the curve above. 

$\cA$-discriminants are notoriously large in all but a few restricted
settings. For instance, the polynomial $\belta_{\{0,404,405,808\}}$ defining 
the curve above has the following coefficient for 
$c^{808}_2c_4$:\linebreak 
\mbox{}\hfill
903947086576700909448402875044761267196347419431440828445529608410806270
\hfill\mbox{} \\
\mbox{}\hfill ...[2062 digits omitted]...
...93441588472666704061962310429908170311749217550336.\hfill \mbox{}\\
Fortunately, we have the following theorem, describing a one-line
parametrization of $\nabla_\cA$.
\begin{hk}
{\em \scalebox{.96}[1]{(See \cite{kapranov}, \cite{passare}, and 
\cite[Prop.\ 4.1]{dfs}.)}\linebreak 
\scalebox{1}[1]{Given $\cA\!:=\!\{a_1,\ldots, a_m\}\!\subset\!\Zn$ with 
$\nabla_\cA$ a hypersurface, the discriminant locus
$\nabla_\cA$ is exactly}\linebreak 
the closure of\\
\scalebox{.95}[1]
{$\left\{ \left[u_1\lambda^{a_1}:\cdots:u_m \lambda^{a_m}\right]\; \left| \; 
u\!:=\!(u_1,\ldots,u_m)\!\in\!\C^m,\
\cA u \!=\!\bO, \ \sum^m_{i=1}u_i\!=\!0, \
\lambda\!=\!(\lambda_1,\ldots,\lambda_n)\!\in\!\Csn\right. \right\}$. \qed}}
\end{hk}

\noindent
Thus, once we know the null-space of an $(n+1)\times m$ matrix, 
we have a formula parametrizing $\nabla_\cA$. 
Recall that for any two subsets $U,V\subseteq\!\R^N$, their {\em Minkowski 
sum} $U+V$ is\linebreak $\{u+v\; | \; u\!\in\!U \ , \ v\!\in\!V\}$.   
Also, for any matrix $M$, we let $M^T$ denote its transpose. 
\begin{cor} 
\label{cor:mink} 
Following the notation above, let $\hA$ denote the 
$(n+1)\times m$ matrix whose $i\thth$ column has coordinates corresponding
to $1\times a_i$, let $B\!\in\!\R^{m\times p}$ be any real 
matrix whose columns are a basis for the right null-space of $\hA$, and 
define $\varphi : \C^{p} \longrightarrow  
\R^{m}$ via $\varphi(t)\!:=\!\log\left|tB^T\right|$. Then 
$\amoeba(\Delta_\cA)$ is the Minkowski sum 
of the row space of $\hA$ and 
$\varphi\!\left(\C^{p}\right)$. \qed 
\end{cor} 

\noindent 
If one is familiar with elimination theory, then it is evident from the 
Horn-Kapranov Uniformization that discriminant amoebae are subspace bundles 
over a lower-dimensional amoeba. This is a geometric reformulation of the 
homogeneities satisfied by the polynomial $\Delta_\cA$. 
\begin{ex}
\label{ex:tetracell}
{\em Continuing Example~\ref{ex:tetra2}, we observe that 
$\hA\!=\!\renewcommand\arraystretch{1.0}\begin{bmatrix}
1 & 1 & 1 & 1 \\ 0 & 404 & 405 & 808 \end{bmatrix}$ has right null-space
generated by the respective transposes of $(1,-405,404,0)$ and $(1,-2,0,1)$. 
The Horn-Kapranov Uniformization then tells us that
$\nabla_\cA$ is simply the closure of the rational surface 
$\left.\left\{ \left[ t_1+t_2 : 
(-405t_1-2t_2) \lambda^{404} : 
404t_1 \lambda^{405} : t_2 \lambda^{808} \right] \; \right| \; 
t_1,t_2\!\in\!\C , \lambda\!\in\!\Cs\right\}\subset 
\Pro^3_\C$. Note that $f$ and $\frac{1}{c_1}f$ have the same roots
and that $u\mapsto u^{1/405}$ is a well-defined bijection on
$\R$ that preserves sign. Note also that the roots of $f$ and
$\barf(y)\!:=\!\frac{1}{c_1}f\!\left(\text{\scalebox{1}[.7]
{$\left(\frac{c_1}{c_3}\right)^{1/405}$}}y\right)$ differ only by a scaling 
when $f$ has real coefficients, and that
$\barf$ is of the form $1+c'_2y^{404}+y^{405}+c'_4y^{808}$.
It then becomes clear that we can reduce the
study of $\nabla_\cA(\R)$ to a lower-dimensional slice: 
intersecting $\nabla_\cA$ with the plane defined by
$c_1\!=\!c_3\!=\!1$ yields the following parametrized curve in 
$\C^2$:\linebreak  
$\babla_\cA\!:=\!\left\{ 
\left. \left(\frac{-405t_1-2t_2}{t_1+t_2} 
\left(\frac{404t_1}{t_1+t_2}\right)^{-404/405}, \ \  
\frac{t_2}{t_1+t_2}
\left(\frac{404t_1}{t_1+t_2}\right)^{-808/405} \right) \; 
\right| \; t_1,t_2\!\in\!\C\right\}$.  
In other words, the preceding curve is the closure of the set of all 
$(c'_2,c'_4)\!\in\!(\Cs)^2$ such that  
$1+c'_2x^{404}+x^{405}+c'_4x^{808}$ has a degenerate 
root in $\Cs$. Our preceding illustration of the image of 
$\babla_\cA(\R)$ within
$\amoeba(\belta_\cA)$ (after taking log absolute values of coordinates) 
thus has the following explicit parametrization:\\
\scalebox{.835}[1]{$\left\{ 
\left(\log|405t_1+2t_2|-\frac{1}{405}\log|t_1+t_2|-
\frac{404}{405}\log|404t_1|, 
\log|t_2|+\frac{403}{405}\log|t_1+t_2|-\frac{808}{405} 
\log|404t_1|\right)\right\}_{[t_1:t_2]\in\Pro^1_\R}$} \dia} 
\end{ex} 

A geometric fact about amoebae that will prove quite useful here is the 
following elegant quantitative result of Passare and Rullg\aa{}rd. Recall 
that the {\em Newton polytope} of a Laurent polynomial 
$f\!\in\!\C[x^{\pm 1}_1,\ldots,x^{\pm 1}_n]$ is the {\em convex 
hull} of\footnote{i.e., smallest convex set containing...}  
the exponent vectors appearing in the monomial term expansion of $f$. 
\begin{prt} \cite[Cor.\ 1]{passarerullgard} 
{\em Suppose $f\!\in\!\C[x^{\pm 1}_1,x^{\pm 1}_2]$ has Newton polygon $P$. 
Then $\area(\amoeba(f))\!\leq\!\pi^2\area(P)$. \qed } 
\end{prt} 

\subsection{Discriminant Chambers and Cones} 
$\cA$-discriminants are central in real root counting because 
the real part of $\nabla_\cA$ determines where,  
in coefficient space, the real zero set of a polynomial changes topology. 
Recall that a {\em cone} in $\R^m$ is any subset closed under 
addition and nonnegative linear combinations. The dimension of a 
cone $C$ is the dimension of the smallest flat containing $C$. 

\begin{dfn} 
\label{dfn:reduction} {\em  
Suppose $\cA\!=\!\{a_1,\ldots,a_m\}\!\subset\!\Zn$ has cardinality $m$ and 
$\nabla_\cA$ is a hypersurface. We then call any connected component $\cC$ 
of the complement of $\nabla_\cA$ in $\Pro^{m-1}_\R\setminus\{c_1\cdots 
c_m\!=\!0\}$ a {\em (real) discriminant chamber}. 
Also let $\hA$ denote the $(n+1)\times m$ matrix whose $i\thth$ column has 
coordinates corresponding to $1\times a_i$ and let $B$ be any real matrix 
$B\!=\![b_{i,j}]\!\in\!\R^{m\times p}$ with $\begin{bmatrix}
\hA\\B^T\end{bmatrix}$ invertible. If $\log|\cC|B$ contains 
an $m$-dimensional cone then we call $\cC$ an {\em outer} chamber 
(of $\nabla_\cA$). All other chambers of $\nabla_\cA$ are called {\em inner} 
chambers (of $\nabla_\cA$). Finally, 
we call the formal expression $(c_1,\ldots,c_m)^B\!:=\!\left(c^{b_{1,1}}_1
\cdots c^{b_{m,1}}_m, \ldots,c^{b_{1,p}}_1\cdots c^{b_{m,p}}_m\right)$ a 
monomial change of variables, and we refer to images of the form $\cC^B$ (with 
$\cC$ an inner or outer chamber) as {\em reduced} chambers. \dia}  
\end{dfn}

\noindent 
It is easily verified that $\log|\cC^B|\!=\!\log|\cC|B$. The latter 
notation simply means the image of $\log|\cC|$ under right multiplication 
by the matrix $B$. 

\begin{ex} {\em 
The illustration from Example \ref{ex:tetra2} shows $\R^2$ 
partitioned into what appear to be $3$ convex and unbounded regions, and 
$1$ non-convex unbounded region. There are in fact $4$ convex and unbounded 
regions, with the fourth visible only if the downward pointed spike were 
allowed to extend much farther down. So $\cA\!=\!\{0,404,405,808\}$ results 
in exactly $4$ reduced outer chambers. \dia } 
\end{ex} 

Note that exponentiating by any $B$ as above yields a well-defined 
multiplicative homomorphism from $(\Rs)^m$ to $(\Rs)^p$ when $B$ has 
rational entries with all denominators odd. In particular, 
the definition of outer chamber is independent of $B$, since (for the 
$B$ considered above) $\Log\!\left|\cC^B\right|$ is unbounded and convex iff 
$\Log\!\left|\cC^{B^*}\right|$ is 
unbounded and convex, where $B^*$ is any matrix whose columns are a basis for 
the orthogonal complement of the row space of $\hA$. 

One can in fact reduce the study of the topology of sparse polynomial 
real zero sets to representatives coming from reduced discriminant chambers. 
A special case of this reduction is the following result. 
\begin{lemma} 
(See \cite[Prop.\ 2.17]{drrs}.)  
\label{lemma:cell} 
Suppose $\cA\!\subset\!\Zn$ has $n+3$ elements, $\cA$ is not contained in 
any $(n-1)$-flat, $\cA\cap Q$ has cardinality 
$n$ for all facets $Q$ of $\conv \cA$, 
all the entries of $B\!\in\!\Q^{(n+3)\times 2}$ have odd denominator, 
and $\begin{bmatrix}\hA\\ B^T\end{bmatrix}$ is invertible.  
Also let $f,g\!\in\!\cF_\cA\setminus\nabla_\cA$ have respective real 
coefficient vectors $c$ and $c'$ with $c^B$ and $c'^B$ lying in the same 
reduced discriminant chamber. Then  
all the complex roots of $f$ and $g$ are non-singular, and 
the respective zero sets of $f$ and $g$ in $(\Rs)^n$ are diffeotopic. 
In particular, for $n\!=\!1$, we have that either $f(x)$ and $g(x)$ 
have the same number of positive roots or $f(x)$ and $g(-x)$ have the 
same number of positive roots. \qed
\end{lemma}

\subsection{Integer Linear Algebra and Linear Forms in Logarithms}  
\label{sub:bit} 
Certain quantitative results on integer matrix factorizations and 
linear forms in logarithms will prove crucial for our main algorithmic 
results. 
\begin{dfn}
\label{dfn:smith} 
{\em 
Let $\Z^{n\times m}$ denote the set of $n\times m$ matrices
with all entries integral, and let $\gln(\Z)$ denote the
set of all matrices in $\Z^{n\times n}$ with determinant $\pm 1$
(the set of {\em unimodular} matrices).
Recall that any $n\times m$ matrix $[u_{ij}]$ with
$u_{ij}\!=\!0$ for all $i\!>\!j$ is called {\em upper triangular}.}

{\em 
Given any $M\!\in\!\Z^{n\times m}$, we then call an
identity of the form $UM=H$, with $H\!=\![h_{ij}]\!\in\!\Z^{n\times m}$
upper triangular and $U\!\in\!\gln(\Z)$, a {\em Hermite factorization}
of $M$. Also, if we have the following conditions in addition:
\begin{enumerate}
\item{$h_{ij}\!\geq\!0$ for all $i,j$.}
\item{\scalebox{.96}[1]{for all $i$, if $j$ is the smallest $j'$ such that
$h_{ij'}\!\neq\!0$ then $h_{ij}\!>\!h_{i'j}$ for all $i'\!\leq\!i$.}}
\end{enumerate}
then we call $H$ {\em \underline{the} Hermite normal form} of $M$. \dia} 
\end{dfn}
\begin{prop}
\label{prop:mono}
We have that $x^{AB}\!=\!(x^A)^B$ for any $A,B\!\in\!\Z^{n\times n}$.
Also, for any field $K$, the map defined by $m(x)\!=\!x^U$,
for any unimodular matrix $U\!\in\!\Z^{n\times n}$, is an automorphism of
$(K^*)^n$. \qed  
\end{prop}
\begin{thm}
\label{thm:unimod}
\cite[Ch.\ 6, Table 6.2, pg.\ 94]{storjophd}
For any $A\!=\![a_{i,j}]\!\in\!\Z^{n\times m}$ with $m\!\geq\!n$, the Hermite 
factorization of $A$ can be computed within $O\!\left(nm^{2.376} 
\log^2(m\max_{i,j}|a_{i,j}|)\right)$ bit operations.
Furthermore, the entries of all matrices in the Hermite 
factorization have bit size $O(m\log(m\max_{i,j}|a_{i,j}|))$. \qed 
\end{thm}

The following result is a very special case of work of Nesterenko that 
dramatically refines Baker's famous theorem on linear 
forms in logarithms \cite{baker}. 
\begin{thm} 
\label{thm:log} 
(See \cite[Thm.\ 2.1, Pg.\ 55]{nesterenko}.)  
For any integers $\gamma_1,\alpha_1,\ldots,\gamma_N,\alpha_N$ with
$\alpha_i\!\geq\!2$ for all $i$, define $\Lambda(\gamma,\alpha)\!:=\!
\gamma_1\log(\alpha_1)+\cdots+\gamma_N\log(\alpha_N)$. Then
$\Lambda(\gamma,\alpha)\!\neq\!0\Longrightarrow 
\log\left|\frac{1}{\Lambda(\gamma,\alpha)}\right|$ is bounded above by 
$2.9(N+2)^{9/2}(2e)^{2N+6}
(2+\log\max_j|\gamma_j|)\prod\limits^N_{j=1}\log \alpha_j$. \qed 
\end{thm} 

The most obvious consequence of lower bounds for linear forms in 
logarithms is an efficient way to determine the signs of 
monomials in integers. 
\begin{algor}\label{algor:sign}\mbox{}\\
{\bf Input:} Positive integers $\alpha_1,u_1,\ldots,\alpha_M,u_M$
and $\beta_1,v_1,\ldots,\beta_N,v_N$ with $\alpha_i,\beta_i\!\geq\!2$ for all
$i$.\\
{\bf Output:} The sign of $\alpha^{u_1}_1\cdots\alpha^{u_M}_M -
\beta^{v_1}_1\cdots\beta^{v_N}_N$. \\
{\bf Description:}\\
\vspace{-.5cm}
\begin{enumerate}
\addtocounter{enumi}{-1}
\item{Check via gcd-free bases (see, e.g., \cite[Sec.\ 8.4]{bs}) whether
$\alpha^{u_1}_1\cdots\alpha^{u_M}_M\!=\!\beta^{v_1}_1\cdots\beta^{v_N}_N$. 
If so, output ``{\tt They are equal.}'' and {\tt STOP}.}
\item{Let $U\!:=\!\max\{u_1,\ldots,u_M,v_1,\ldots,v_N\}$ and\\ 
\mbox{}\hfill $E\!:=\!\frac{2.9}{\log 2}(2e)^{2M+2N+6}(1+\log U)\times
\left(\prod\limits^M_{i=1}\log \alpha_i 
\right)\left(\prod\limits^N_{i=1}\log \beta_i \right)$.\hfill\mbox{} }
\item{For all $i\!\in\![M]$ (resp.\ $i\!\in\![N]$), let $A_i$ (resp.\ $B_i$)
be a rational number agreeing with $\log\alpha_i$ (resp.\ $\log \beta_i$) in
its first $2+E+\log_2M$ (resp.\ $2+E+\log_2N$) leading bits.\footnote{For
definiteness, let us use Arithmetic-Geometric Mean Iteration
as in \cite{dan} to find these approximations.}}
\item{Output the sign of
$\left(\sum\limits^M_{i=1}u_iA_i\right)-\left(\sum\limits^N_{i=1}v_iB_i
\right)$ and {\tt STOP}.}
\end{enumerate}
\end{algor}

\begin{lemma}
\label{lemma:bin}
Algorithm \ref{algor:sign} is correct and, following the preceding notation, 
runs within a number of bit operations asymptotically linear in\\
\mbox{}\hfill $(M+N)(30)^{M+N}L(\log U)
\left(\prod\limits^M_{i=1}L(\log \alpha_i)\right)
\left(\prod\limits^N_{i=1}L(\log \beta_i)\right)$,\hfill\mbox{}\\
where $L(x)\!:=\!x(\log x)^2\log\log x$.  \qed
\end{lemma}

\noindent
Lemma \ref{lemma:bin} follows immediately from 
Theorem \ref{thm:log}, the well-known fast iterations for 
approximating log (see, e.g., \cite{dan}), and the known refined bit 
complexity estimates for fast multiplication (see, e.g., 
\cite[Table 3.1, pg.\ 43]{bs}).

\section{Chamber Cones and Polyhedral Models} 
\label{sec:chambercones} 
\subsection{Defining and Describing Chamber Cones} 
\begin{dfn} {\em Suppose $X\!\subset\!\R^m$ is convex and 
$Q\!\supseteq\!X$ is the polyhedral cone consisting of all 
$c\!\in\!\R^m$ with $c+X\!\subseteq\!X$. We call 
$Q$ the {\em recession cone} for $X$ and, if $p\!\in\!\R^m$ 
satisfies\linebreak (1) $p+Q\!\supseteq\!X$ and (2) $p+c+Q\!\not\supseteq\!X$ 
for any $c\!\in\!Q\setminus\{\bO\}$, then we call $p+Q$ the 
{\em placed} recession cone. In particular, 
the placed recession cone for $\Log|\cC|$ for $\cC$ an outer  
chamber (resp.\ a reduced outer chamber) is called 
a {\em chamber cone} (resp.\ a {\em reduced chamber cone}) of $\nabla_\cA$. 
We call the facets of the (reduced) chamber cones of $\nabla_\cA$ {\em 
(reduced) walls of $\nabla_\cA$}. We also refer to walls of dimension $1$ as 
{\em rays}. \dia} 
\end{dfn}

\noindent 
\begin{minipage}[b]{0.55\linewidth}
\vspace{0pt}
\begin{ex} 
{\em
Returning to Example \ref{ex:tetra2}, let us 
draw the rays that are the boundaries of the $4$  
reduced chamber cones: While there appear to be 
just $3$ reduced chamber cones, there are in fact $4$,  
since there is an additional (slender) reduced 
chamber cone with vertex placed much farther down. 
(The magnified illustration to the right actually shows 
$2$ nearly parallel rays going downward, very close together.) 
Note also that reduced chamber cones need not share vertices. \dia } 
\end{ex} 
\end{minipage}\hspace{.2cm}
\begin{minipage}[b]{0.35\linewidth}
\vspace{0pt}
\scalebox{2}[1]{\epsfig{file=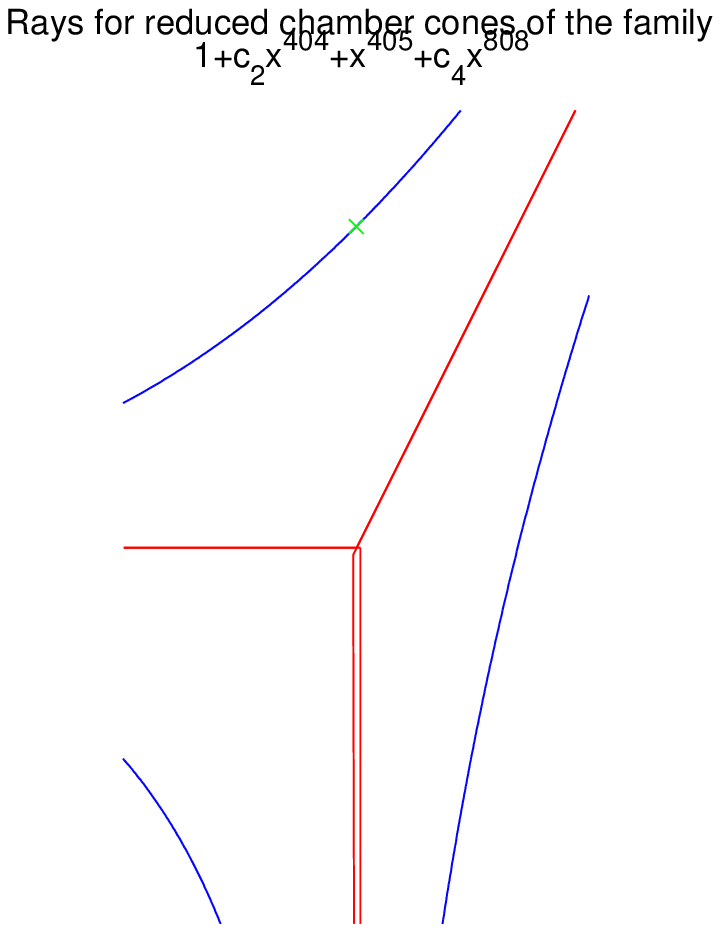,height=1.83in,clip=}} 
\end{minipage} 

\noindent
Chamber cones are well-defined since chambers are {\em log-convex}, 
being the domains of convergence of a particular class of hypergeometric 
series (see, e.g., \cite[Ch.\ 6]{gkz94}). A useful corollary 
of the Horn-Kapranov Uniformization is a surprisingly simple and explicit 
description for chamber cones. 
\begin{dfn}
\label{dfn:radiant} 
{\em  
Suppose $\cA\!\subset\!\Z^n$ has cardinality $n+3$, is not 
contained in any $(n-1)$-flat, and is not a pyramid.\footnote{The last 
assumption simply means that 
there is no point $a\!\in\!\cA$ such that $\cA\setminus\{a\}$ lies in an 
$(n-1)$-flat.} 
Also let $B$ be any real $(n+3)\times 2$ matrix whose columns are a
basis for the right null space of $\hA$ and let $\beta_1,\ldots,\beta_{n+3}$ 
be the rows of $B$. We then call any set of indices $\cI\!\subset\!\{1,
\ldots,n+3\}$ satisfying: \\
\mbox{}\hspace{1.5cm}(a) $[\beta_i]_{i\in \cI}$ is a maximal 
rank $1$ submatrix of $B$.\\ 
\mbox{}\hspace{1.5cm}(b) $\sum_{i\in\cI} \beta_i$ is 
{\em not} the zero vector.\\ 
a {\em radiant subset corresponding to $\cA$}.  \dia } 
\end{dfn} 
\begin{thm}
\label{thm:walls}
Suppose $\cA\!\subset\!\Z^n$ has cardinality $n+3$, is not 
contained in any $(n-1)$-flat, is not a pyramid, and $\nabla_\cA$ is a 
hypersurface. 
Also let $B$ be any real $(n+3)\times 2$ matrix whose columns are a
basis for the right null space of $\hA$ and let $\beta_1,\ldots,\beta_{n+3}$ 
be the rows of $B$. Finally, let 
$S\!:=\![s_{i,j}]\!:=\!B\begin{bmatrix}0 & -1\\ 1 & 0
\end{bmatrix} B^T$ and let $s_i$ denote the row vector whose 
$j\thth$ coordinate is $0$ or $\log|s_{i,j}|$ according as 
$s_{i,j}$ is $0$ or not. Then each wall of 
$\nabla_\cA$ is the Minkowski sum of the 
row-space of $\hat{\cA}$ and a ray of the form $s_i-\R_+\sum_{j\in \cI} 
e_j$ for some unique radiant subset $\cI$ of $\cA$ and any $i\!\in\!\cI$. 
In particular, the number of walls of $\nabla_\cA$, the number of chamber 
cones of $\nabla_\cA$, and the number of radiant subsets corresponding to 
$\cA$ are all identical, and at most $n+3$.  
\end{thm} 

\noindent 
Note that the hypotheses on $\cA$ are trivially satisfied when 
$\cA\!\subset\!\Z$ has cardinality $4$. Note also that the definition 
of a radiant subset corresponding to $\cA$ is independent of the 
chosen basis $B$ since the definition is invariant under column 
operations on $B$. 
\begin{rem} 
{\em Theorem \ref{thm:walls} thus refines an earlier result of 
Dickenstein, Feichtner, and Sturmfels (\cite[Thm.\ 1.2]{dfs}) 
where, in essence, {\em un}shifted variants of chamber cones (all going 
through the origin) were computed for nonpyramidal $\cA\!\subset\!\Zn$ of 
arbitrary cardinality. A version of Theorem \ref{thm:walls} for $\cA$ of 
arbitrary cardinality will appear in \cite{prrt}. \dia } 
\end{rem} 
\begin{ex} {\em 
It is easy to show that a generic $\cA$ satisfying the hypotheses of 
our theorem will have exactly $n+3$ chamber cones, e.g., 
Example \ref{ex:tetra2}. It is also 
almost as easy to construct examples having fewer chamber cones. 
For instance, taking $\cA\!:=\!\left\{\begin{bmatrix} 0\\0\end{bmatrix},
\begin{bmatrix} 1\\0\end{bmatrix},
\begin{bmatrix} 2\\0\end{bmatrix},
\begin{bmatrix} 0\\1\end{bmatrix},
\begin{bmatrix} 2\\1\end{bmatrix}\right\}$, we see that 
$B\!:=\!\text{\scalebox{.7}[.5]{$\begin{bmatrix}[r] -1 & 0\\ 2 & 1\\ -2 & -2\\ 
0 & 1\\ 1 & 0\end{bmatrix}$}}$ satisfies the hypotheses of Theorem 
\ref{thm:walls} and that $\{1,5\}$ is a {\em non}-radiant subset. So the 
underlying $\nabla_\cA$ has only $3$ chamber cones. \dia}  
\end{ex} 

\noindent
{\bf Proof of Theorem \ref{thm:walls}:} 
First note that by Corollary \ref{cor:mink}, $\amoeba(\Delta_\cA)$ is the 
Minkowski sum of $\varphi\!\left(\C^2\right)$ and 
the row space of $\hA$, where $\varphi(t)\!:=\!\Log\left|tB^T\right|$. 
So then, determining the walls reduces to determining the directions 
orthogonal to the row space of $\hA$ in which $\varphi(t)$ becomes unbounded. 

Since $\mathbbm{1}\!:=\!(1,\ldots,1)$ 
is in the row space of $\hA$, we have $\mathbbm{1}B\!=\!\bO$ and thus 
$\varphi(t)\!=\!\varphi(t/M)$ for all $M\!>\!0$. So we can 
restrict to the compact subset $\{(t_1,t_2)\; | \; 
|t_1|^2+|t_2|^2\!=\!1\}$ and observe that $\varphi(t)$ becomes unbounded iff 
$t\beta^T_i\longrightarrow 0$ for some $i$. In particular, we see that there 
are no more than $n+3$ reduced walls. Note also that 
$t\beta^T_i\longrightarrow 0$ 
iff $t$ tends to a suitable (nonzero) multiple of $\beta_i\begin{bmatrix} 
0 & -1 \\ 1 & 0\end{bmatrix}$, in which case the coordinates of 
$\varphi(t)$ becoming unbounded are exactly those with index $j\!\in\cI$ 
where $\cI$ is the unique radiant subset corresponding to those rows of 
$B$ that are nonzero multiples of $b_i$. (The assumption that 
$\cA$ not be a pyramid implies that $B$ can have no zero rows.) 
Furthermore, the coordinates of $\varphi(t)$ that become unbounded 
each tend to $-\infty$. Note that Condition (b) of the radiance 
condition comes into play since we are looking for directions 
{\em orthogonal to the row-space of $\hA$} in which 
$\varphi(t)$ becomes unbounded. 

We thus obtain that each wall is of the 
asserted form. However, we still need to account for the coordinates of 
$\varphi(t)$ that remain bounded. Now, if $t$ tends to a suitable (nonzero) 
multiple of $\beta_i\begin{bmatrix} 0 & -1 \\ 1 & 0\end{bmatrix}$, then 
it clear that any coordinate of $\varphi(t)$ of index $j\!\not\in\!\cI$ 
tends to $s_{i,j}$ (modulo a multiple of $\mathbbm{1}$ added to 
$\varphi(t)$). So we have indeed described every wall, and given 
a bijection between radiant subsets corresponding to $\cA$ and 
the walls of $\nabla_\cA$. 

To conclude, note that the row space of $\hA$ has dimension $n+1$ by 
construction, so the 
walls are all actually (parallel) $n$-plane bundles over rays. So by the 
Archimedean Amoeba Theorem, each outer chamber of $\nabla_\cA$ must be 
bounded by $2$ walls and the walls have a natural cyclic ordering. 
Thus the number of chamber cones is the same as the number of rays and we are 
done. \qed 

\subsection{Which Chamber Cone Contains Your Problem?} 
An important consequence of Theorem \ref{thm:walls} is that while the
underlying $\cA$-discriminant polynomial $\Delta_\cA$ may have
huge coefficients, the {\em rays} of a linear projection of 
$\amoeba(\Delta_\cA)$ admit a concise description involving few bits, save 
for the transcendental coordinates coming from the ``shifts'' $s_i$. 
By applying our quantitative estimates from 
Section \ref{sub:bit} we can then quickly find which chamber cone contains 
a given $n$-variate $(n+3)$-nomial.  
\begin{thm}
\label{thm:lp}
Following the notation of Theorem \ref{thm:walls}, suppose $f\!\in\cF_\cA
\cap \R[x_1,\ldots,x_n]$ and let $\tau$ denote the maximum 
bit size of any coordinate of $\cA$. Then we can determine the unique chamber 
cone containing $f$ --- or correctly decide if $f$ is contained in $2$ or more 
chamber cones --- within a number of arithmetic operations polynomial in  
$n+\tau$. Furthermore, if $f\!\in\cF_\cA \cap \Z[x_1,\ldots,x_n]$, 
$\sigma$ is the maximum bit size of any coefficient of $f$,  
and $n$ is fixed, then we can also obtain a {\em bit} complexity bound 
polynomial in $\tau+\sigma$. \qed 
\end{thm}

Theorem \ref{thm:lp} is the central tool behind our complexity results 
and follows immediately upon proving the correctness of (and giving 
suitable complexity bounds for) the following algorithm:   
\begin{algor}\label{algor:conesort}\mbox{}\\
{\bf Input:} A subset $\cA\!\subset\!\Zn$ of cardinality $n+3$ (with 
$\cA$ not a pyramid and not contained in
any\linebreak
\mbox{}\hspace{1.5cm}$(n-1)$-flat, and $\nabla_\cA$ a hypersurface) and the 
coefficient vector $c$ of an\linebreak 
\mbox{}\hspace{1.5cm}$f\!\in\!\cF_\cA\cap \R[x_1,\ldots,x_n]$.\\ 
{\bf Output:} Radiant subsets $\cI$ and $\cI'$ (corresponding to 
$\cA$) generating the walls of the unique\linebreak
\mbox{}\hspace{1.85cm}chamber cone containing 
$f$, or a true declaration that $f$ is contained in at least $2$\linebreak  
\mbox{}\hspace{1.85cm}chamber cones.\\  
{\bf Description:}\\
\vspace{-.5cm}
\begin{enumerate}
\addtocounter{enumi}{-6}
\item{(Preprocessing) Compute the Hermite Factorization 
$H^T\!=\!U^T\hA$ and let $B$ be the \linebreak submatrix defined by the 
rightmost $2$ columns of $U$.} 
\item{(Preprocessing) 
Let $\beta_1,\ldots,\beta_{n+3}$ be the rows of $B$, 
$S\!:=\![s_{i,j}]\!:=\!B\begin{bmatrix}0 & -1\\ 1 & 0
\end{bmatrix} B^T$, and let $s_i$ denote the row vector whose
$j\thth$ coordinate is $0$ or $\log|s_{i,j}|$ according as
$s_{i,j}$ is $0$ or not.} 
\item{(Preprocessing) Find all radiant subsets $\cI\!\subset\!\{1,\ldots,
n+3\}$ corresponding to $\cA$.}
\item{(Preprocessing) For any radiant subset $\cI$ let 
$\beta'_\cI\!:=\!-\sum_{j\in \cI}\beta_j$ and let 
$s_\cI$ denote the row vector $s_iB$ for any fixed $i\!\in\!\cI$.} 
\item{(Preprocessing) Sort the $\beta'_\cI$ in order 
of increasing counter-clockwise angle with the $x$-coordinate ray and 
let $R$ denote the resulting ordered collection of $\beta'_\cI$.}  
\item{(Preprocessing) For any radiant subset $\cI$ let 
$v_\cI\!\in\!\Q^2$ denote the intersection of the lines 
$s_\cI+\R\beta'_\cI$ and $s_{\cI'}+\R\beta'_{\cI'}$, where 
$\beta'_{\cI'}$ is the counter-clockwise neighbor of $\beta'_{\cI}$.} 
\item{Let $ConeCount\!:=\!0$.} 
\item{Via binary search, attempt to find a pair of adjacent rays of the 
form\\ 
$(v_\cI+\R_+\beta'_\cI,v_\cI+\R_+\beta'_{\cI'})$ containing $\Log|c|B$.
\begin{enumerate}
\item{If ($ConeCount\!=\!0$ and there is no such pair of rays)\\ 
or\\ 
($ConeCount\!=\!1$ and there is such a pair of rays)\\ 
then output 
{\tt ``Your $f$ lies in at least $2$ distinct chamber cones.''} and 
{\tt STOP}.} 
\item{\scalebox{.94}[1]{If $ConeCount\!=\!0$ and there is such a pair of 
rays, set $ConeCount\!:=\!ConeCount+1$,}\linebreak  
delete $\beta'_\cI$ and $\beta'_{\cI'}$ from $R$, and {\tt GOTO STEP} (2).} 
\end{enumerate}} 
\item{Output {\tt ``Your $f$ lies in the unique chamber cone determined by 
$\cI$ and $\cI'$.''} and {\tt STOP}.} 
\end{enumerate} 
\end{algor} 
\begin{rem} 
{\em An important detail for large scale computation is that the 
preprocessing steps (-5)--(0) need only be done {\em once} 
per support $\cA$. This can significantly increase efficiency in applications 
where one has just one (or a few) $\cA$ and one needs to answer chamber 
cone membership queries for numerous $f$ with the same support. \dia } 
\end{rem}

\noindent 
{\bf Proof of Correctness of Algorithm \ref{algor:conesort}:} First note 
that the computed matrix $B$ indeed has columns that form a basis for 
the right null-space of $\cA$. This is because our assumptions on 
$\cA$ ensure that the rank of $\hA$ is $n+1$ and thus the last $2$ 
rows of $H^T$ consist solely of zeroes. 

By construction, Theorem \ref{thm:walls} then implies that 
the $\beta'_\cI$ are exactly the reduced rays for $\nabla_\cA$, modulo 
an invertible linear map. (The invertible map arises because 
right-multiplication by $B$ induces an injective but non-orthogonal 
projection of the right null-space of $\hA$ onto $\R^2$.) 

It is then clear that the preprocessing steps do nothing more than 
give us a $B$ suitable for Theorem \ref{thm:walls} and a sorted 
set of reduced rays ready for chamber cone membership 
queries via binary search. In particular, since the reduced chamber cones 
cover $\R^2$, the correctness of Steps (1)--(3) is clear and we are done. 
\qed 

\medskip 
\noindent 
In what follows, we will use the ``soft-Oh'' notation $O^*(f)$ to 
abbreviate bounds of the form $O\!\left(f(\log f)^{O(1)}\right)$. 

\medskip 
\noindent 
{\bf Complexity Analysis of Algorithm \ref{algor:conesort}:} 
Let us first approach our analysis from the more involved point of 
view of bit complexity. Our arithmetic complexity bound will then 
follow upon a quick revisit. 

By Theorem \ref{thm:unimod} (which quotes bounds from \cite{storjophd}), 
Step (-5) takes $O\!\left(n^{3.376}\tau^2\right)$ bit operations. Also, the 
resulting bit size for the entries of $B$ is $O(n\tau)$. 

The complexity of Step (-4) is negligible, save for the 
approximation of certain logarithms. The latter won't come into play until 
we start checking on which side of a ray a point lies, so let us 
analyze the remaining preprocessing steps. 

Step (-3) can be done easily through a greedy approach: one simply 
goes through the rows $\beta_2,\ldots,\beta_{n+3}$ to see which rows are 
multiples of $\beta_1$. Once this is done, one checks if the resulting 
set of indices is indeed radiant or not, and then one repeats this process 
with the remaining rows of $B$. In summary, this entails $O(n^2)$ arithmetic 
operations on number of bit-size $O(n\tau)$, yielding a total of 
$O^*(n^3\tau)$ bit operations. 

Step (-2) has negligible complexity.

The comparisons in Step (-1) can be accomplished by computing the cosine and 
sine of the necessary angles via dot products and cross products. Via 
the well-known asymptotically optimal sorting algorithms, it is then 
clear that Step (-1) requires $O(n\log n)$ arithmetic operations on 
integers of bit size $O(n\tau)$, contributing a total of 
$O^*(n^2\tau\log n)$ bit operations. 

Step (0) has negligible complexity. 

At this point, we see that the complexity of the Preprocessing Steps 
(-5)--(0) is $O(n^{3.376}\tau^2)$ bit operations.  

Continuing on to Steps (1)--(3), we now see that we are faced with 
$O(\log n)$ sidedness comparisons between a point and an 
oriented line. More precisely, we need to evaluate $O(\log n)$ signs of 
determinants of matrices of the form $\begin{bmatrix}\Log|c|B-s_\cI\\ 
\beta'_\cI \end{bmatrix}$. Each such sign evaluation, thanks to Algorithm 
\ref{algor:sign} and Lemma \ref{lemma:bin}, takes 
$O\!\left(n30^{2n+5}L(\sigma+n\tau)L(\sigma)^{n+3}L(n\tau)^{n+2}\right)$ 
bit operations. 

We have thus proved our desired bit complexity bound which, while 
polynomial in $\tau+\sigma$ for fixed $n$, is visibly 
exponential in $n$. Note however that the exponential bottleneck occurs 
only in the sidedness comparisons of Step (2). 

To obtain an improved arithmetic complexity 
bound, we then simply observe that the sidedness comparisons can be 
replaced by computations of signs of differences of monomials, 
simply by exponentiating the resulting linear forms in logarithms. 
Via recursive squaring \cite[Thm.\ 5.4.1, pg.\ 103]{bs}, it is then clear 
that each such comparison requires only $O(n^2\tau)$ arithmetic operations. 
So the overall number of arithmetic operations drops to polynomial 
in $n+\tau$ and we are done.  \qed 

\medskip 
Let us now state some final combinatorial constructions before fully 
describing how chamber cones apply to real root counting. 

\subsection{Canonical Viro Diagrams and the Probability of 
Lying in Outer Chambers}
\label{sub:can} 
Our use of outer chambers and chamber cones enables us to augment an earlier 
construction 
of Viro. Let us first recall that a {\em triangulation} of a
point set $\cA$ is simply a simplicial complex $\Sigma$ whose vertices lie in
$\cA$. 
\begin{dfn}{\em 
We say that a triangulation of $\cA$ is {\em coherent} iff
its maximal simplices are exactly the domains of linearity for some
function $\ell$ that is convex, continuous, and piecewise linear on the convex 
hull of $\cA$. In 
particular, we will sometimes define such an $\ell$ by fixing 
the values $\ell(a)$ for just those $a\!\in\!\cA$ and then employing the 
convex hull of the points $(a,\ell(a))$ as $a$ ranges over $\cA$. The resulting 
graph is known as the {\em lower hull} of the {\em lifted} point set 
$\{(a,\ell(a))\; | \; a\!\in\!\cA\}$. \dia } 
\end{dfn} 
\begin{dfn}
\label{dfn:viro} {\em  
(See Proposition 5.2 and Theorem 5.6 of
\cite[Ch.\ 5, pp.\ 378--393]{gkz94}.)
Suppose $\cA\!\subset\!\Zn$ is finite and the convex hull of $\cA$ has
positive volume and boundary $\partial A$. Suppose also that $\cA$
is equipped with a coherent triangulation $\Sigma$ and a
function $s : \cA \longrightarrow \{\pm\}$ which we will call a
{\em distribution of signs for $\cA$}. We then call any  
edge with vertices of opposite sign an {\em alternating edge}, and we define a
piece-wise linear manifold  --- the {\em Viro diagram}
$\cV_\cA(\Sigma,s)$ --- in the following local manner: For any $n$-cell
$C\!\in\!\Sigma$,
let $L_C$ be the convex hull of the set of all midpoints of alternating edges of
$C$, and then define $\cV_\cA(\Sigma,s)\!:=\!\bigcup\limits_{C \text{ an } 
n\text{-cell}} L_C\setminus\partial A$. When $\cA\!=\!\supp(f)$ and $s$ is 
the corresponding sequence of coefficient signs, then we also call
$\cV_{\Sigma}(f)\!:=\!\cV_\cA(\Sigma,s)$ the {\em Viro diagram of $f$ 
corresponding to $\Sigma$}. \dia } 
\end{dfn}
\begin{ex} 
\label{ex:penta} 
{\em Consider 
$f(x)\!:=\!1-x_1-x_2+3x^4_1x_2+3x_1x^4_2$. 
Then
$\supp(f)\!=\!\{(0,0),(1,0),$\linebreak
$(0,1),(1,4),(4,1)\}$ and has convex hull a pentagon.
So then there are exactly $5$ coherent triangulations, yielding $5$ possible
Viro diagrams for $f$ (drawn in thicker green lines):\\
\mbox{}\hfill\epsfig{file=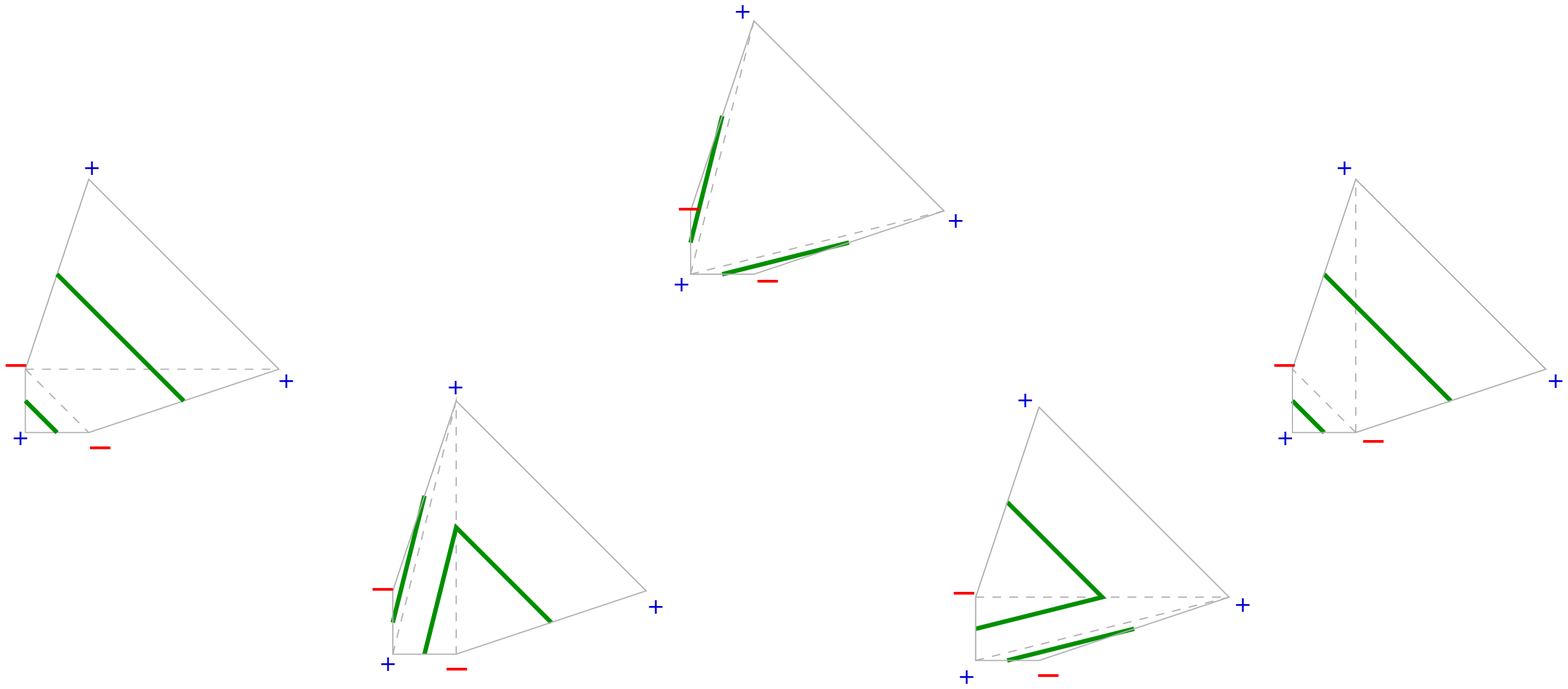,height=2.5in}  \hfill \mbox{} \\
Note that all these diagrams have exactly $2$ connected components, with
each component isotopic to an open interval. Note also that our $f$ here is a
$2$-variate $(2+3)$-nomial. \dia} 
\end{ex}

\begin{dfn} 
\label{dfn:lift} 
{\em 
Suppose $\cA\!\subset\!\Z^n$ has cardinality $n+3$,
is not a pyramid, is not contained in any $(n-1)$-flat, and $\nabla_\cA$ 
is a hypersurface.
Also let $B$ be any real $(n+3)\times 2$ matrix whose columns are a
basis for the right null space of $\hA$. For any $f\!\in\!\cF_\cA$ let us then 
define $v(f)\!:=(v_1(f),\ldots,v_{n+3}(f))\!:=\!\left(\sum_{i\in\cI} 
e_i\right)+ \left(\sum_{j\in\cJ} e_j\right)$ where $\cI$ and $\cJ$ are the 
unique radiant subsets corresponding to the unique chamber cone containing 
$\Log|c|$. (We set $v(f)\!:=\!\bO$ should there not be a unique such 
chamber cone.) Let us then define $\nanewt(f)$ to be the convex hull of 
$\{(a_i,v_i)\; | \; i\!\in\!\{1,\ldots,n+3\}\}$ and let 
$\Sigma_f$ denote the triangulation of $\cA$ induced by the lower 
hull of $\nanewt(f)$.  
Also, we call any polynomial of the form $\sum\limits_{a_i\in Q} 
c_ix^{a_i}$ --- with $Q$ a cell of $\Sigma_f$ ---  a {\em canonical lower 
polynomial} for $f$. Finally, we call $\cV(f)\!:=\!\cV_\Sigma(f)$ the 
{\em canonical Viro diagram} of $f$. \dia } 
\end{dfn} 
\begin{ex} 
{\em Let $f(x_1)\!:=\!1-\frac{1}{2}x^{404}+x^{405}-2x^{808}$, 
$c\!:=\!\left(1,-\frac{1}{2},1,2\right)$, and 
$\cA\!:=\!\{0,404,405,808\}$. A routine calculation then reveals that 
$\{\{2\},\{3\}\}$ is the pair of radiant subsets corresponding to 
the unique chamber cone containing $\Log|c|$. We then obtain that 
$v(f)\!=\!(0,1,1,0)$ and thus $\nanewt(f)$ is \epsfig{file=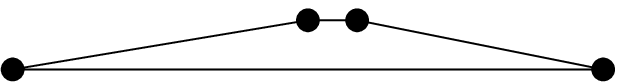,
height=.15in} (modulo some artistic stretching). In particular, $\Sigma_f$ 
thus has the single cell 
$[0,808]$, which is an alternating cell, and thus $\cV(f)$ consists of a 
single point. More than coincidentally, $f$ has exactly $1$ positive 
root. \dia } 
\end{ex} 
\begin{ex} 
{\em Returning to Example \ref{ex:penta}, let $c\!:=\!(1,-1,-1,
3,3)$. A routine calculation then reveals that the unique chamber cone 
containing $\Log|c|$ is defined by the pair of radiant subsets 
$\{\{2\},\{3\}\}$. So then $v(f)\!=\!(0,1,1,0,0)$ and $\Sigma_f$ 
is the exactly the upper middle triangulation from the illustration of 
Example \ref{ex:penta}. $\cV(f)$ then consists of $2$ disjoint 
open intervals and, more than coincidentally, the positive zero set 
of $f$ has exactly $2$ connected components, each homeomorphic to an 
open interval. \dia } 
\end{ex} 
\begin{thm}
\label{thm:vtx}
Following the notation above, set $\hat{f}_t(x)\!:=\!\sum^{n+3}_{i=1}c_i
t^{v_i(f)} x^{a_i}$ and assume in addition that 
$\cA\cap Q$ has cardinality $n$ for all facets $Q$ of $\conv \cA$. 
Then $c\!:=\!(c_1,\ldots,c_{n+3})$ lies in an outer 
chamber $\Longrightarrow$ the positive zero sets of $\hat{f}_t$, 
as $t$ ranges over $(0,1]$, are each diffeotopic to the positive 
zero set of $\hat{f}_1$. In particular, $\hat{f}_1\!=\!f$ and thus, when $c$ 
lies in an outer chamber, the positive zero set of 
$f$ is isotopic to $\cV(f)$. 
\end{thm}
\begin{rem}
{\em 
For $n\!=\!1$ we thus obtain that the number of positive 
roots of the tetranomial $f$ is exactly the cardinality of its 
canonical Viro diagram. \dia }
\end{rem} 

\noindent
{\bf Proof:} By construction, the image of 
$\Log|(c_1t^{v_1(f)},\ldots,c_{n+3}t^{v_{n+3}(f)})|$ as $t$ ranges over 
$(0,1]$ is simply a ray entirely contained in a unique chamber cone. Moreover, 
by assumption (and since outer chambers are log convex), the ray is also 
contained entirely in $\Log|\cdot|$ of an outer chamber. So the 
first part of our theorem follows from Lemma \ref{lemma:cell}. 

The final part of our theorem is then just  
a reformulation of Viro's Theorem on 
the isotopy type of toric deformations of real algebraic sets 
(see, e.g., \cite[Thm.\ 5.6]{gkz94}). \qed 

\medskip 
The main contribution of our paper is thus an efficient method to associate 
a {\em canonical} Viro diagram to the positive zero set of a given $f$ so that 
both $C^1$ manifolds have the same topology. Such a method appears to be 
new, although the necessary ingredients have existed in the literature 
since at least the 1990s. In particular, to the best of our knowledge, 
all earlier applications of Viro's method 
designed clever $f$ having the same topology as some specially 
tailored Viro diagram, 
thus going in the opposite direction of our construction.  

We state up front that our method does {\em not} work for all $f$. 
However, our development yields a sufficient condition 
--- outer chamber membership --- that holds with high 
probability under the stable log-uniform measure. 
\begin{thm} 
\label{thm:nplus3} Suppose $\cA\!\subset\!\Z^n$ has cardinality $n+3$,
is not a pyramid, is not contained in any $(n-1)$-flat, and $\nabla_\cA$
is a hypersurface. Suppose also that the coefficients of an 
$f\!\in\!\cF_\cA\cap \R[x_1,\ldots,x_n]$ are independently chosen 
via the stable log-uniform measure over $\R$ (resp.\ $\Z$). 
Then with probability $1$, $f$ lies in some outer chamber. 
In particular, if we assume in addition that $\cA\cap Q$ has cardinality $n$ 
for all facets $Q$ of $\conv \cA$, then the positive zero set of 
$f$ is isotopic to $\cV(f)$ with probability $1$. 
\end{thm} 

\noindent 
{\bf Proof:} By Theorem \ref{thm:walls}, $\amoeba(\Delta_\cA)$ 
is an $n$-plane bundle over $\amoeba\!\left(\belta_\cA\right)$, where 
$\belta_\cA\!\in\!\Z[a,b]$ and 
$\Delta_\cA(c_1,\ldots,c_{n+3})\!=\!\belta_\cA(\alpha(c),\beta(c))$ 
for suitable monomials $\alpha$ and $\beta$ in $c$. Furthermore, 
thanks to Corollary 8 of \cite{sadykovsolid}, the 
complement of $\amoeba(\belta_\cA)$ has no bounded convex connected 
components. 

Letting $c$ denote the coefficient vector of $f$ we thus obtain that 
$f$ lies in an outer chamber iff $\Log|c|\!\not\in\!\amoeba(\Delta_\cA)$. 
In particular, it is clear by the Passare-Rullg\aa{}rd Theorem 
that in any large centered cube $C$, 
the volume of $\amoeba(\Delta_\cA)\cap C$ occupies a vanishingly  
small fraction of $C$. So the first assertion is proved. 

The final assertion is an immediate consequence of the first assertion 
and Theorem \ref{thm:vtx}. \qed  

Theorem \ref{thm:pos} then follows easily from Theorems \ref{thm:vtx} and 
\ref{thm:nplus3}. 
The applications of Theorems \ref{thm:vtx} and \ref{thm:nplus3} to 
computational real topology will be pursued in another paper. 

\subsection{Proving Theorem \ref{thm:pos}} 
\label{sub:log} 
Let us first consider the following algorithm for counting the 
positive roots of ``most'' real univariate tetranomial. 
\begin{algor}\label{algor:tetra}\mbox{}\\
{\bf Input:} A tetranomial $f\!\in\!\R[x_1]$ with support 
$\cA$.\\
{\bf Output:} A number in $\{0,1,2,3\}$ that, for any $f$ in an 
outer chamber of $\nabla_\cA$, is exactly the number of positive roots of 
$f$.\\
{\bf Description:}\\
\vspace{-.5cm}
\begin{enumerate}
\item{Via Algorithm \ref{algor:conesort}, and any sub-quadratic 
planar convex hull algorithm (see, e.g., \cite{compugeom}), compute the 
canonical Viro diagram $\cV(f)$. }  
\item{If $f$ did not lie in a unique chamber cone then 
output {\tt ``Your $f$ does not lie in an outer chamber so please 
use an alternative method.''} and {\tt STOP}. } 
\item{Output the cardinality of $\cV(f)$ and {\tt STOP}. } 
\end{enumerate} 
\end{algor} 

Assuming Algorithm \ref{algor:tetra} is correct, we can count the real 
roots of $f$ simply by applying Algorithm \ref{algor:tetra} to $f(x_1)$ 
and $f(-x_1)$. Theorem \ref{thm:pos} thus follows immediately upon proving 
the correctness of our last algorithm and a suitable complexity bound. 
We now do this. 

\medskip 
\noindent 
{\bf Proof of Correctness of Algorithm \ref{algor:tetra}:} 
By Theorem \ref{thm:vtx}, the number of 
positive roots of $f$ is exactly the cardinality of $\cV(f)$ 
whenever $f$ is in an outer chamber. So we indeed have correctness. \qed 

\medskip 
\noindent 
{\bf Complexity Analysis of Algorithm \ref{algor:tetra}:} 
Let us first observe that Algorithm \ref{algor:tetra} indeed 
gives a correct answer with probability $1$ (relative to the 
stable log-uniform measure): this is immediate from Theorem 
\ref{thm:nplus3}. 

So now we need only prove a suitable complexity bound. Let us 
begin with the more refined setting of bit complexity: From our 
earlier complexity analysis of Algorithm \ref{algor:conesort}, it is clear 
that Step (1) needs $O(\log^2 D)+
O\!\left(L(h+\log D)L(h)^{4}L(\log D)^{3}\right)$ bit operations, 
neglecting the computation of $\cV(f)$. The complexity of 
computing $\cV(f)$ (which is essentially dominated by computing the 
convex hull of $4$ points with coordinates of bit size $O(\log D)$) 
is clearly negligible in comparison. So we obtain a final 
bit complexity bound of $O^*\!\left((h+\log D)h^4\log^3 D\right)$. 
(The complexity of Steps (2) and (3) is negligible.) 

For arithmetic complexity, our earlier complexity analysis of Algorithm 
\ref{algor:conesort} specializes easily to an upper bound of 
$O(\log^2 D)$. (The speed-up arises from the easiness of 
checking inequalities involving integral powers of real numbers 
in the BSS model over $\R$.) So we are done. \qed 

\begin{rem} {\em 
It is important to note that when $f$ lies in a chamber cone but 
{\em not} in any outer chamber, Algorithm \ref{algor:tetra} can 
give a wrong answer. However, thanks to Theorem \ref{thm:nplus3}, 
such an occurence has probability $0$ under the stable log-uniform 
measure. \dia}  
\end{rem}  

\section{Proving Theorem \ref{thm:neg}}
To prepare for the proof, we first set up some notation.  Fix positive 
integers $\ell$ and $m$.  Let $P = (p_{ij}) \in \N^{\ell \times m}$ be a matrix of natural numbers ordered as \[ p_{11} \geq p_{21} \geq \cdots \geq p_{s1} \ \ \text{ and } \ \  p_{ij} > p_{i(j+1)},\]
for $i = 1,\ldots,s$ and $j = 1, \ldots, m-1$.  Also, let $a_{ij}$ be 
indeterminates over the same indexing set.  Consider now the following 
polynomial:
\begin{equation}\label{supposedsos}
\begin{split}
S_P(x) := \ &  \sum_{i=1}^{\ell} 
{\left(\sum_{j=1}^m  a_{ij}x^{p_{ij}}\right)^2}  \in \N[a_{ij}][x] \\
= \ & g_d(P) x^d + \cdots + g_1(P) x + g_0(P), \\
\end{split}
\end{equation}
in which each nonzero $g_i(P)$ is a homogeneous (quadratic) polynomial in 
$\N[a_{ij}]$.  
Notice that at most $\ell m^2$ of these $g_i$ are nonzero.  We will refer to 
the integer $p_{ij}$ as the \textit{exponent} corresponding to the 
\textit{coefficient} $a_{ij}$.

\begin{lemma}\label{finiteglemma}
The set of polynomials $g_i(P)$ arising from $S_P(x)$ as $P$ ranges over all 
$\ell \times m$ nonnegative integer matrices is finite.
\end{lemma}

\noindent 
{\bf Proof:} 
The form of Equation (\ref{supposedsos}) allows for at most $2^{lm^2}$ 
possible such $g_i$. \qed 

Suppose now that $f \in \mathbb R[x]$ has degree $d$ and is a sum of $\ell$ 
squares, each involving at most $m$ terms.  Then, there is a set of 
exponents $P$ and an assignment $\bar{a}_{ij} \in \R$ for the coefficients 
$a_{ij}$ such that  $f  = S_P(x)$.  If we fix a set of exponents $P$ for $f$, 
then conversely, any real point in the variety determined by the $g_i$ gives a 
representation of $f$ as a sum of $\ell$ squares, each involving at most $m$ 
terms.

We will prove Theorem \ref{thm:neg} using contradiction by showing that a 
certain infinite family of trinomials cannot all have sparse representations 
of the form in Equality (\ref{supposedsos}).  For this approach to work, 
however, we will 
need to find a universal set of coefficients $\bar{a}_{ij}$ that ``represents" 
all of them.

\begin{lemma}\label{universalaij}
Let $F \subset \R[x]$ be an infinite collection of polynomials which are 
sums of $\ell$ squares, each involving at most $m$ terms.  Moreover, suppose 
that the nonzero coefficients of polynomials $f \in F$ come from a finite set 
$C$.  Then, there is an infinite subset $f_1,f_2,\ldots \in F$ of them with 
corresponding exponents $P_1,P_2,\ldots \in \N^{\ell \times m}$ and a single 
set of real coefficients $\bar{a}_{ij}$  such that $f_k = S_{P_k}(x) |_{a_{ij} 
= \bar{a}_{ij}}$ for all $k$.  
\end{lemma}

\noindent 
{\bf Proof:} 
Let $f_1,f_2,\ldots \in F$ be an infinite sequence of polynomials with 
corresponding exponents $P_1,P_2,\ldots \in \N^{\ell \times m}$.  Also let 
$S$ be the set of all possible (polynomial) coefficients of polynomials 
$S_{P_k}(x)$; from Lemma \ref{finiteglemma}, this set is finite.  By 
assumption, each $f_k$ gives rise to a set of equations involving a subset of 
$S$ and real numbers coming from the finite set $C$.  The set of all such 
equations is again finite, and therefore, by the infinite pigeon-hole 
principle, there is a subsequence of the $f_k$ which have the same set of 
equations governing their coefficients $a_{ij}$.  Picking any real solution to 
such a set of equations finishes the proof.

We now have all the tools we need to prove Theorem \ref{thm:neg}.  However, to make the argument less cumbersome, we first recall ``little-oh" notation: for a 
function $h: \N \to \R$, we say that $h(n)\!=\!o(n)$ 
if \[\lim_{n \to \infty} \frac{h(n)}{n} = 0.\]
It is easy to see that the sum of any finite number of such functions is also 
$o(n)$.  Moreover, if $\lim_{n \to \infty} \frac{p(n)}{n} = p$ for some 
constant $p$, then $p(n) = np + o(n)$. 

\medskip 
\noindent 
{\bf Proof of Theorem \ref{thm:neg}:} 
Suppose, by way of contradiction, that every positive definite trinomial can be written sparsely as the sum of $\ell$ squares, each involving at most $m$ terms.  Consider the infinite sequence of positive definite trinomials, 
\begin{equation}\label{trinomialfamily}
f_k = x^{2k} + x^{2k-1} + 1, \ \ \ k = 1,2,\ldots.
\end{equation}
Using Lemma \ref{universalaij}, we can find a subsequence $f_{k_s}$ with 
corresponding exponents $P_{k_1},P_{k_2},\ldots \in \N^{\ell \times m}$ and a 
single set of real numbers $\bar{a}_{ij}$ such that $f_{k_s} = S_{P_{k_s}}(x) 
|_{a_{ij} = \bar{a}_{ij}}$ for $s = 1,2,\ldots$.  We will assume that the 
$\bar{a}_{ij}$ are chosen such that the number of them that are zero is 
maximal among all such sets of coefficients $\bar{a}_{ij}$.  Moreover, we 
shall assume that this maximality still holds if we restrict to any infinite 
subsequence of the $f_{k_s}$.
For clarity of exposition in what follows, we will not keep updating the  
subscripting of indices when taking subsequences.

Given an exponent matrix $P_{k_s} \in \N^{\ell \times m}$, define a new matrix \[\tilde{P}_{k_s} = \frac{1}{k_s} P_{k_s}.\] 
This corresponds naturally to the transformation $x \mapsto x^{1/{k_s}}$ applied to both sides of an equation $f_{k_s}(x) = S_{P_{k_s}}(x)$.
Since $\deg(f_{k_s}) = 2k_s$, each matrix $\tilde{P}_{k_s}$ has entries in the 
interval $[0,1]$.  By compactness, we may choose a subsequence $P_{k_s}$ such 
that $\tilde{P}_{k_s}$ converges in the (entry-wise) Euclidean norm to a matrix $\tilde{P} = (\tilde{p}_{ij}) \in [0,1]^{\ell \times m}$.  Clearly, we have $\tilde{p}_{11} = 1$ and also that some entry of $\tilde{P}$ is $0$; it turns out 
that these are the only possible numbers for entries of $\tilde{P}$ which play 
a role in (\ref{supposedsos}).

\textit{Claim}: If  $\tilde{p}_{ij} \neq 0,1$, then $\bar{a}_{ij} = 0$.

Suppose, on the contrary, that $\tilde{P}$ contains $r > 2$ real numbers, 
$p_1,\ldots,p_r$ with $0= p_1 < \cdots < p_r = 1$.
Each power of $x$ occurring in a summand of Equality (\ref{supposedsos}),  
after squaring, is of the form 
\begin{equation}\label{exponentssqring} 
k_s p_u + k_s p_v + o(k_s).
\end{equation}
Thus, for all sufficiently large $s$, the powers of $x$ occurring in 
Equality (\ref{supposedsos}) are partitioned into classes determined by the 
number of distinct values, \[p_u + p_v, \ \  u,v = 1,\ldots,r.\]  

Notice that the numbers in Formula (\ref{exponentssqring}) all become strictly 
smaller 
(resp.\ larger) than \mbox{$2k_s-1$} (resp. $0$) as $s \to \infty$ unless $
u =v = r$ (resp. $u=v=1$).  In particular, with $s$ large, the polynomials 
$g_{w}(P_{k_s}) \in \N[a_{ij}]$ for 
\begin{equation}\label{beginendgs}
w = 2k_s + o(k_s) \text{ \ and \ } w = 0 + o(k_s)
\end{equation}
do not involve the indeterminates $a_{ij}$ coming from exponents of the form 
$k_s p_u + o(k_s)$ with $u \neq 1,r$. Conversely, each monomial in every 
polynomial $g_{w}(P_{k_s})$ with $w$ not in one of  the classes from 
(\ref{beginendgs}) 
contains at least one such indeterminate $a_{ij}$.  

Since the only nonzero coefficients of the sequence (\ref{trinomialfamily}) 
come from the classes from (\ref{beginendgs}),
it follows that we may replace with $0$ all coefficients $\bar{a}_{ij}$ 
corresponding to exponents $k_s p_u + o(k_s)$ with $u \neq 1,r$ and still 
have an equality \[f_{k_s} = S_{P_{k_s}}(x)|_{a_{ij} = \bar{a}_{ij}}.\]  The 
claim therefore follows from the maximality property of the chosen set of 
coefficients $\bar{a}_{ij}$.

We next  examine what all this says about the limiting behavior of the 
expressions from Equality (\ref{supposedsos}).  From the claim, it follows 
that when $s$ is 
large, we need only consider those exponents from the matrices $P_{k_s}$ that 
are on the order
\[ k_s + o(k_s) \text{ \ and \ } 0 + o(k_s).\]  Fix such a large $s$ and 
consider the smallest exponent $p$ on the order $k_s + o(k_s)$ that occurs 
with a nonzero coefficient in (\ref{supposedsos}) after substituting the 
$\bar{a}_{ij}$ for the $a_{ij}$.  When the sum from (\ref{supposedsos}) is 
expanded, the term $x^{2p}$ will appear with positive  coefficient; i.e., 
\[g_{2p}(P_{k_s}) |_{a_{ij} = \bar{a}_{ij}} \neq 0.\]  It follows from the 
form of the sequence (\ref{trinomialfamily}) that $p = k_s$ (since the other 
nonzero term is odd).  In particular, it is not possible to obtain a nonzero 
coefficient for $x^{2k_s-1}$ in $f_{k_s}$.  This contradiction completes the 
proof. \qed 

\bibliographystyle{amsalpha}

\end{document}